\documentclass[11pt]{article}
\usepackage{amsmath,amssymb}

\newtheorem{Definition}{Definition}

\newtheorem{Lemma}{Lemma}

\newcommand{\tmop}[1]{\ensuremath{\operatorname{#1}}}

\newtheorem{Theorem}{Theorem}

\newtheorem{Problem}{Problem}
\newtheorem{Conjecture}{Conjecture}
\newenvironment{proof}{
  \noindent\textbf{Proof}\ }{\hspace*{\fill}
  \begin{math}\Box\end{math}\medskip}

\begin{document}

\title{Quadrance polygons, association schemes and strongly regular graphs}
\author{Le Anh Vinh\\School of Mathematics\\University of New South Wales\\Email: vinh@maths.unsw.edu.au} \maketitle

\begin{abstract}
	Quadrance between two points $A_1 = [x_1,y_1]$ and $A_2 = [x_2,y_2]$ is the number $Q ( A_1, A_2 ) := ( x_2 - x_1 )^2 + ( y_2 - y_1 )^2$. In this paper, we present some interesting results arise from this notation. 
In Section 1, we will study geometry over finite fields under quadrance notations. The main purpose of Section 1 is to answer the question, for which $a_1,\ldots a_n$, we have a polygon $A_1\ldots A_n$ such that $Q(A_i,A_{i+1})=a_i$ for $i = 1,\ldots,n$. In Section 2, using tools developed in Section 1, we define a family of association schemes over finite field space $F_q \times F_q$ where $q$ is a prime power. These schemes give rise to a graph $V_q$ with vertices the points of $F_q^2$, and where $(X,Y)$ is an edge of $V_q$ if and only if $Q(X,Y)$ is a nonzero square number in $F_q$. In Section 3, we will show that $V_q$ is a strongly regular graph and propose a conjecture about the maximal clique of $V_q$.

\end{abstract}

\section{Universal geometry over finite fields}

Suppose that $q$ is an odd prime power, and that $F_q$ is the finite field
with $q$ elements. To avoid lengthly calculations (but mostly repetition), we assume that $q$ is of form $q = 4l+3$ for some integer $l$ throughout Sections 1-3. In Section \ref{companion results}, we will present companion results for case $q = 4l + 1$. The following definitions follow \cite{quadrance} where the importance of the notation of quadrance is developed.

\begin{Definition}
  The \textbf{quadrance} $Q ( A_1, A_2 )$ between the points $A_1 = [ x_1,
  y_1]$, and $A_2 = [ x_2, y_2 ]$ is the number
  \[ Q ( A_1, A_2 ) := ( x_2 - x_1 )^2 + ( y_2 - y_1 )^2 . \]
\end{Definition}

\begin{Definition}
  A \textbf{circle} in a finite field $F_q$ with center $A_0 \in F_q \times F_q$ and quadrance $K \in F_q$ is set of all points $X$in $F_q \times F_q$ such that
  \[ Q(A_0, X) = K. \]
\end{Definition}

Recall that a (multiplicative) character of $F_q$ is a homomorphism from
$F^{\ast}$, the multiplicative group of the non-zero elements of $F_q$, to the
multiplicative group of complex numbers with modulus $1$. The identically $1$
function is the principal character of $F_q$ and is denoted $\chi_0$. Since
$x^{q - 1} = 1$ for every $x \in F_q^{\ast}$ we have $\chi^{q - 1} = \chi_0$
for every character $\chi$. A character $\chi$ is of order $d$ if $\chi^d =
\chi_0$ and $d$ is the smallest positive integer with this property. By
convention, we extend a character $\chi$ to the whole of $F_q$ by putting
$\chi ( 0 ) = 0$. The quadratic (residue) character is defined by $\chi ( x )
= x^{( q - 1 ) / 2}$. Equivalently, $\chi$ is $1$ on square, $0$ at $0$ and $-
1$ otherwise. It is easy to see that $\sum_{i \in F_q} \chi ( i ) = 0$. Let
$A_q, B_q, C_q$ and $D_q$ be the numbers of $i$ in $F_q$ such that $\chi ( i )
= \chi ( i + 1 ) = 1$, $\chi ( i ) = \chi ( i + 1 ) = - 1$, $\chi ( i ) = -
\chi ( i + 1 ) = 1$ and $\chi ( i ) = - \chi ( i + 1 ) = - 1$, respectively. We have the following lemma.

\begin{Lemma} \label{squarefinite} Let $q$ be an odd prime power of the form $q = 4l + 3$ for some integer $l$. Then
$A_q = B_q = (q-3)/4$, $C_q = (q+1)/4$ and $D_q = (q-3)/4$. 
\end{Lemma}

\begin{proof}
We have
  \begin{align*}
    \sum_{i \in F_q} \left( 1 + \chi(i) \right) \left( 1
    + \chi(i+1) \right) &= q + \sum_{i \in F_q} \left\{  \chi(i) + \chi(i+1) 
    + \chi(i ( i + 1 )) \right\}\\
    &= q + \sum_{i \in F_q, \neq 0, - 1} \chi(i^2 ( 1 + i^{- 1} ))\\
    &= q + \sum_{i \in F_q, \neq 0, - 1} \chi(1 + i^{- 1})\\
    &= q - \chi(1) - \chi(0)
    \\   
    &= q - 1. 
  \end{align*}
But we have
  \begin{equation*}
   \left( 1 + \chi(i) \right) \left( 1 + \chi(i + 1)\right) =
   \begin{cases}
    4  & \text{if both}\; \; i \;\text{and}\; i+1 \; \text{are square}  ,\\
    2  & \text{if}\; \; i = 0.\\
    0 & \text{otherwise}.
   \end{cases}
  \end{equation*}
  Thus, $4 A_q + 2 = q - 1$ or $A_q = ( q - 3 ) / 4$.  
  Similarly, consider the sum
  \begin{equation*}
  \sum_{i \in F_q} \left( 1 - \chi(i) \right) \left(
     1 - \chi(i + 1) \right). 
  \end{equation*}
  This gives $4 B_q + 1 = q - 2$ or $B_q = ( q - 3 ) / 4$.
  Similarly, consider the sum
  \begin{equation*}
  \sum_{i \in F_q} \left( 1 + \chi(i) \right) \left(
     1 - \chi(i + 1) \right), 
  \end{equation*}
  we have $C_q = (q+1)/4$. And consider the sum   
  \begin{equation*}
  \sum_{i \in F_q} \left( 1 - \chi(i) \right) \left(
     1 + \chi(i + 1) \right), 
  \end{equation*}
  we have $D_q = (q-3)/4$. This concludes the proof.
\end{proof}

From Lemma \ref{squarefinite}, we can count the number of points in any circle in
$F_q^2$.

\begin{Lemma} \label{number points}
  Let $C$ be a circle with nonzero quadrance in $F_q^2$. Then $C$ has exactly $N = q+1$ points.
\end{Lemma}
\begin{proof} Without lose of generality, we may assume that $C$ is centered at the origin. Let $K$ be the quadrance of $C$. We will prove the lemma for the case $K = 1$ as other cases are similar.
  The number of points in $C$ is
  \[ N = \left| \left\{ [ x, y ] \in F_q^2 \mid x^2 + y^2 = 1 \right\} \right|. \]
  We have $x^2 + y^2 = 1$ if and only if $x^2 = ( - y^2 ) + 1$. Since $- y^2$ is not a square in $F_q^2$.
    We have the numbers of $( x^2, y^2 )$ with $x, y \neq 0$ is $D_q = ( q - 3 )
    / 4$. For each $( x^2, y^2 )$ with $x, y \neq 0$, we have $4$
    corresponding points $[x, y ]$. Besides, we have 4 other points $[ 0, 1 ], [ 0, - 1 ],
    [ - 1, 0 ]$ and $[ 1, 0 ]$. Hence, $N = 4 D_q + 4 = q + 1$.    
  This concludes the proof of the lemma.
\end{proof}

Let $C_i(X)$ denote the circle centered at $X \in F_q^2$ with quadrance $i$.

\begin{Lemma} \label{khac khong}
Let $i, j \ne 0$ in $F_q$ and let $X, Y$ be two distinct points in $F_q^2$ such that $k = Q(X,Y)\ne 0 $. Then $\left|C_i(X) \cap C_j(Y)\right|$ only depends on $i, j$ and $k$. Precisely,  define \[f ( i, j, k ) := i j - \frac{ (i - j - k )^2 }{ 4}.\]
Then the number of intersection points is $p_{i j}^k$, where
  \begin{equation}\label{51}
  p_{ij}^k =
   \begin{cases}
    0  & \text{if}\; \; f(i,j,k)\; \text{is non-square},\\
    1  & \text{if}\; \; f(i,j,k) = 0,\\
    2 &\text{if}\; \; f(i,j,k)\; \text{is square}.
   \end{cases}
  \end{equation}

\end{Lemma}

\begin{proof}
Suppose that $X = [ m, n ]$ and $Y = [ m + x, n + y ]$ for some $m, n, x, y \in F_q$ then $x^2 + y^2 = k$. Suppose that $Z \in C_i(X) \cap C_j(Y)$ where $Z = [ m + x + u, n + y + v ]$ for some $u, v \in F_q$. Then we have $u^2 + v^2 = j$ and $( x + u )^2 + ( y + v^2 ) = i$. This implies that $x u + y v = (i - j - k)/2$.
  But we have $( x u + y v )^2 + ( x v - y u )^2 = ( x^2 + y^2 ) ( u^2 + v^2
  )$ so
  \[ ( x v - y u )^2 = k j - \frac{( i - j - k )^2}{4} = i j - \frac{( k - i - j )^2}{4} = f ( i, j, k ) . \]
  If $f ( i, j, k )$ is non-square number in $F_q$ then it is clear that there does not exist
  such $x, y, u, v$, or $p_{i j}^k = 0$. Otherwise, let $\alpha = ( i - j
  - k ) / 2$ and $f ( i, j, k ) = \beta^2$ for $0 \leqslant \beta \leqslant (
  p + 1 ) / 2$ then
  \[
    x v - y u = \pm \beta,\qquad x u + y v = \alpha.
  \]
  Solving for $( u, v )$ with respect to $( x, y )$ we have
  \[
    u = ( \alpha x \mp \beta y ) / k,\qquad  v = ( \alpha y \pm \beta x) / k. 
  \]
  If $\beta = 0$ then we have only one $( u, v )$ for each $( x, y )$, but if
  $\beta \neq 0$ then we have two pairs $( u, v )$. This implies (\ref{51}), completing the proof.
\end{proof}

The function $f(i,j,k)$ defined above has an interesting property.

\begin{Lemma} \label{function}
Given $i, j \ne 0$ in $F_q$, there exists at least $|C_i|/2$ values of $k$ such that $f(i,j,k)$ is square (maybe zero).
\end{Lemma}
\begin{proof}
There are two cases.
  \begin{enumerate}
    \item Suppose that $i j$ is square in $F_q$. Then $i j = v^2$ for some $v \in F_q$. We have 
    \[f(i,j,k) = v^2 \left( 1 - \left( \frac{k-i-j}{2v} \right)^2 \right).\]
    From Lemma 1, we have $D_q = ( q - 3 ) / 4$ values of $- t^2$ ($t \neq 0$) such
    that $1 - t^2$ is square.
    For each $t^2 \neq 0$, we have two values of $k = \pm 2 t v + i + j$.
    Besides, we may choose $( k-i-j )$ from $0, \pm 2v$ so we have three more
    values of $k$. Thus, the number of $k$ such that $f ( i, j, k )$ is square
    (maybe zero) is $3 + ( q - 3 ) / 2 = ( q + 3 ) / 2$.
    
    \item Suppose that $i j$ is non square in $F_q$. Then $i j = - v^2$ for some $v \in F_q$. We have 
    \[f(i,j,k) = v^2 \left( - 1 - \left( \frac{k-i-j}{2v} \right)^2 \right).\]
    From Lemma 1, we have $C_q = ( q + 1 ) / 4$ values of $- t^2$ ($t \neq 0$)
    such that $- 1 - t^2$ is square. For each $t^2 \neq 0$, we have two values of $k = \pm 2 t v + i + j$.
    These are all possiblities for $k$. Thus the number of $k$ such that $f ( i, j, k )$ is square (maybe zero) is
    $( q + 1 ) / 2$.
  \end{enumerate}
  This concludes the proof of the lemma.
\end{proof}

The Lemmas \ref{khac khong} and \ref{function} can be applied to yield an interesting result.

\begin{Theorem}\label{polygon}
  If $q \equiv 3 \:( \tmop{mod} 4 )$ then for any $n \geqslant 4$ and $a_1,
  \ldots, a_n \in F_q^{\ast}$, there exists a polygon $A_1 \ldots A_n$ in
  $F_q^2$ such that $Q ( A_i, A_{i + 1} ) = a_i$ for all $i \in \{ 1, \ldots,
  n \}$ (note that $A_{n + 1} = A_1$).
\end{Theorem}

\begin{proof}
  It suffices to prove for $n = 4$. For any $a_1, \ldots, a_4 \in
  F_q^{\ast}$, we need to show that there exists a quadrangle $A_1 A_2 A_3 A_4$
  such that $Q ( A_i, A_{i + 1} ) = a_i$ for $i = 1, \ldots, 4$. From Lemma
  \ref{function}, there are at least $( q + 1 ) / 2$ values of $k$ such that $f ( a_1,
  a_2, k )$ is square and at least $( q + 1 ) / 2$ values of $k$ such that $f
  ( a_3, a_4, k )$ is square. We have only $q$ possible values for $k$ so by
  the pigeonhole principal, there must be some value $k$ such that both $f (
  a_1, a_2, k )$ and $f ( a_3, a_4, k )$ are square. Choose any two points $X,
  Z \in F_q^2$ with $Q ( X, Z )$. From Lemma \ref{khac khong} , $f ( a_1, a_2, k )$ and $f
  ( a_3, a_4, k )$ are square in $F_q$ so there exists $Y, T$ such that $Q (
  X, Y ) = a_1$, $Q ( Y, Z ) = a_2$, $Q ( T, Z ) = a_3$ and $Q ( Z, X ) =
  a_4$. Thus we can choose $A_1 A_2 A_3 A_4 \equiv X Y Z T$. This concludes the proof.
\end{proof}

\section{Quadrance Association Schemes}

We now have enough tools to construct an association scheme on $F_q^2$. Let us recall the formal definition of an association scheme from \cite{association scheme book}.

\begin{Definition}\label{ass scheme}
  An \textbf{association scheme} with $s$ associate classes on a finite
  set $\Omega$ is a partition of $\Omega \times \Omega$ into sets $C_0, C_1,
  \ldots, C_s$ (called \textbf{associate classes}) such that
  \begin{itemize}
    \item[(i)] $C_0 = \tmop{Diag} ( \Omega ) = \{ ( \omega, \omega ) : \omega \in
    \Omega \}$;
    
    \item[(ii)] $C_i$ is symmetric for $i = 1, \ldots, s$, i.e. $C_i' := \{ ( \beta,
    \alpha ) : ( \alpha, \beta ) \in C_i \} = C_i$;
    
    \item[(iii)] for all $i, j, k$ in $\{ 0, \ldots, s \}$ there is an integer $p_{i
    j}^k$ such that, for all $( \alpha, \beta )$ in $C_k$
    \[ | \left\{ \gamma \in \Omega : ( \alpha, \gamma ) \in C_i, ( \gamma,
       \beta ) \in C_j \right\} | = p_{i j}^k . \]
  \end{itemize}
\end{Definition}

The numbers $p_{i j}^k$ are called the \textbf{\textit{intersection numbers}}. The number $a_i = p_{i i}^0$ is called the \textit{\textbf{valency}} of the $i$-th associate class.

Our association schemes are defined by quadrances between elements of $F_q^2$. To meet the condition (i) in Definition \ref{ass scheme}, we need $Q(X, Y) = 0$ if and only if $X \equiv Y$ for all $X, Y \in F_q^2$, which is equivalent to  $q \equiv 3$ (mod $4$). The following theorem gives us an association scheme on $F_q^2$.

\begin{Theorem} \label{quad scheme}
Let $\Omega = F_q^2$, and partition $\Omega \times \Omega$ into
$q$ subsets $C_0, C_1, \ldots, C_{q - 1}$ (indexed by elements of $F_q$; for example $q-1 = -1$ ) such that for $X, Y$ in $F_q^2$, $(
X, Y ) \in C_i$ if and only if $Q ( X, Y ) = i$. This partition gives us an association scheme on $F_q^2$. Furthermore, for $i, j, k \in F_q$, let $f ( i, j, k ) = i j - ( i - j - k )^2 / 4$. Then the
  intersection numbers of the scheme are
  \begin{equation}\label{52}
  p_{ij}^0 =
   \begin{cases}
    0  & \text{if}\; \; i \neq j,\\
    1  & \text{if}\; \; i = j = 0,\\
    q+1 & \text{otherwise},
   \end{cases}
  \end{equation}
and
  \begin{equation}\label{53}
  p_{ij}^k =
   \begin{cases}
    0  & \text{if}\; \; f(i,j,k)\; \text{is non-square},\\
    1  & \text{if}\; \; f(i,j,k) = 0,\\
    2 &\text{if}\; \; f(i,j,k)\; \text{is nonzero square}.
   \end{cases}
  \end{equation}
for $k \neq 0$.
\end{Theorem}

\begin{proof}
  The conditions (i) and (ii) of Definition \ref{ass scheme} are easy to check. We will check
  (iii). Suppose that we have $( X, Y ) \in C_k$ where $X = [ m, n ]$ and $Y = [ m
  + x, n + y ]$ for some $m, n, x, y \in F_p$. Then $x^2 + y^2 = k$. For $k = 0$,
  it is clearly that $x = y = 0$, i.e. $X = Y$. We have
  \begin{equation*}
  \{ Z \in F_p : ( X, Z ) \in C_i, ( Z, Y ) \in C_j \} =
    \begin{cases} 
     \emptyset & \tmop{if} i \neq j,\\
     \{Z \in F_p : Q(X,Z) = i\} & \tmop{if} i = j. 
    \end{cases} 
  \end{equation*}
  This equation and Lemma \ref{number points} imply (\ref{52}). Moreover, Lemma \ref{khac khong} implies (\ref{53}), completing the proof of the theorem.\end{proof}
  
We call this scheme the \textbf{quadrance association scheme}. This scheme can be used to obtain some other schemes. Precisely, we have the following theorem.

\begin{Theorem} \label{induced scheme}
  For any $t$ $|$ $(q - 1)$, there exists an association scheme with $t + 1$
  associate classes on $F_q^2$. 
\end{Theorem}

\begin{proof}
  If $t$ $|$ $q - 1$ then $q - 1 = t n$ for some positive integer $n$. There exists a primitive element of $F_q$, says    $g$. Let $\alpha = g^t$ and let $M$ be set of all
  nonzero $n^{\tmop{th}}$ powers in $F_q$. We have $| M | = t$, and write $M =
  \{ a_1^n, \ldots, a_t^n \}$. Set $x_i = a_i^n$ for $1 \leqslant i \leqslant
  t$. We partition $F_q^2$ into $t + 1$ sets $D_0, D_1, \ldots, D_t$ by
  \[ D_0 = C_0,\qquad D_h = \bigcup_{v = 0}^{n - 1} C_{\alpha^v x_j}\qquad \tmop{for}\qquad h =
     1, \ldots, t. \]
  Here $C_0, C_1, \ldots, C_{q - 1}$ are the associate classes of the quadrance
  association scheme. We check that this partition gives us an association
  scheme. The conditions (i) and (ii) of Definition \ref{ass scheme}
  are easy to check. We will check (iii).
  For $k = 0$, $( X, Y ) \in D_0$ iff $X = Y$. Similarly as in proof of 
  Theorem \ref{quad scheme}, we have that
  \[| \left\{ Z \in F_p : ( X, Z ) \in D_i, ( Z, Y )
  \in D_j \right\} |\] is the same for any $( X, Y ) \in D_0$.
  
  For $k \neq 0$, suppose that $( X, Y ) \in D_k$. Then $( X, Y ) \in
  C_{\alpha^v x_k}$ for some $0 \leqslant v \leqslant n - 1$. The number of
  $Z$'s in $F_q^2$ such that $( X, Z ) \in D_i$ and $( Z, Y ) \in D_j$ is
  \begin{equation} \label{55} V_{i j}^k ( v ) : = \sum_{c = 0}^{n - 1} \sum_{d = 0}^{n - 1} p_{\alpha^c
     x_i, \alpha^d x_j}^{\alpha^v x_k} . \end{equation}
  From Theorem \ref{quad scheme}, $p_{i j}^k$ depends only on whether $4 f ( i, j, k )$ is nonzero square, zero or nonsquare. Thus, we only need
  to look at this condition. For $0 \leqslant v, w \leqslant n - 1$, we have
  \begin{align*}
  \alpha^{2(w-v)}4f(\alpha^c & x_i, \alpha^d x_j, \alpha^v x_k) =	
   \alpha^{2 ( w - v )} \left[ 4 \alpha^{c + d} x_i x_j - ( \alpha^v x_k - \alpha^c x_i - \alpha^d x_j )^2 \right]\\
   &= 4 \alpha^{( c + w - v ) + ( d +
     w - v )} x_i x_j - ( \alpha^w x_k - \alpha^{c + w - v} x_i - \alpha^{d +
     w - v} x_j )^2 ]\\
   &= 4f(\alpha^{(c+w-v)}x_i, \alpha^{(d+w-v)} x_j, \alpha^w x_k).
  \end{align*}
  So that $4f(\alpha^cx_i, \alpha^d x_j, \alpha^v x_k)$ is a square (nonzero/zero) or non square if and only
  if so is $f(\alpha^{(c+w-v)}x_i, \alpha^{(d+w-v)} x_j, \alpha^w x_k)$. Hence, if we replace $( c, d
  )$ by $( c + w - v, d + w - v )$ in (\ref{55}) then $V_{i j}^k ( v ) = V_{i j}^k
  ( w$). This implies that $| \left\{ Z \in F_p : ( X, Z ) \in D_i, ( Z, Y )
  \in D_j \right\} |$ is the same for any $( X, Y ) \in D_k$, completing
  the proof.
\end{proof}

\section{Quadrance Graphs}

\subsection{Strongly Regular Graphs}

Remind that $q \equiv 3\: ( \tmop{mod} 4 )$. In Theorem \ref{induced scheme}, we saw that for any $t$ $|$ $(q - 1)$, there exists an association scheme with $t + 1$ associate classes on $F_q^2$. Now $q$ is an odd prime power so $2$ $|$ $(q - 1)$ and we can compose the quadrance association scheme into an association scheme with three associative classes. It is not hard to see (from the construction in Theorem \ref{induced scheme}) that the association scheme with 3 associative classes is a partition $F_q^2 = \{D_0, D_1, D_2\}$ where 
\begin{align*}
D_0 &= \{(X,X) : X \in F_q\},\\
D_1 &= \{(X,Y) : X, Y \in F_q, Q(X,Y) \tmop{is} \; \tmop{square}\},\\
D_2 &= \{(X,Y) : X, Y \in F_q, Q(X,Y) \tmop{is} \; \tmop{nonsquare}\}.
\end{align*}
Consider a graph of $V_q$ with vertices the points of $F_q^2$, and where $(X,Y)$ is an edge of graph if and only if $(X,Y)$ is in $D_1$. For $X, Y \in F_q^2$, define \[\delta(X,Y) = \{Z : (X,Z), (Z,Y)) \in D_1\}.\]
We have the following theorem.

\begin{Theorem} \label{important} Let $X, Y$ be any two points in $F_q^2$.
\begin{itemize}
\item [(a)] If $(X,Y) \in D_1$ then $|\delta(X,Y)| = (q^2-5)/4$. 
\item [(b)] If $(X,Y) \in D_2$ then $|\delta(X,Y)| = (q^2-1)/4$. 
\end{itemize}
\end{Theorem}
\begin{proof}
  a) Suppose that $( X, Y ) \in D_1$. Then $Q ( X, Y ) = k^2$ for some $k\in F_q^{\ast}$. 
  
  We count the number of points $Z$ in $\delta ( X, Y )$. 
  Let $Z$ be any point in $\delta ( X, Y )$ then there exists $i, j \in F_q^{\ast}$
  such that $Q ( X, Z ) = i^2$ and $Q ( Y, Z ) = j^2$. Without loss of generality, 
  we can divide both $i, j$ by $k$ and assume that $k = 1$. We have
  \begin{align*}
    4 f ( i^2, j^2, 1 ) &= 4 i^2 j^2 - ( 1 - i^2 - j^2 )^2 \\
    &= ( ( i + j )^2 - 1 ) ( 1 - ( i - j )^2 ) . &  & 
  \end{align*}
  There are three separate cases.
  \begin{enumerate}
    \item Suppose that $f ( i^2, j^2, 1 )$ is nonsquare. Then from Lemma \ref{khac khong} , there does
    not exist $Z$ such that $Q ( X, Z ) = i^2$ and $Q ( Y, Z ) = j^2$.
    
    \item Suppose that $f ( i^2, j^2, 1 ) = 0$. We have three subcases.
    \begin{enumerate}
      \item If $( i + j )^2 - 1 = 1 - ( i - j )^2 = 0$ then $i + j = \pm 1$
      and $i - j = \pm 1$. But it contradicts to the condition $i, j \in
      F_q^{\ast}$.
      
      \item If $( i + j )^2 - 1 = 0 \neq 1 - ( i - j )^2$ then $i + j = \pm 1$
      and $( i - j )^2 \neq \pm 1$. Thus, we have $2 ( q - 2 )$ pairs for $( i
      + j, i - j )$. From Lemma \ref{khac khong}, each pair gives us one point $Z$.
      
      \item If $1 - ( i - j )^2 = 0 \neq ( i + j )^2 - 1$ then $i - j = \pm 1$
      and $i + j \neq \pm 1$. Similarly, we have $2 ( q - 2 )$ points $Z$.
    \end{enumerate}
    \item Suppose that $f ( i^2, j^2, 1 )$ is a nonzero square. Then either both $( i + j )^2 - 1$
    and $1 - ( i - j )^2$ are nonzero square or both $( i + j )^2 - 1$ and $1 - ( i - j )^2$ are nonsquare.
       \begin{enumerate}
      \item Suppose that $( i + j )^2 - 1$ and $1 - ( i - j )^2$ are nonzero square. From
      Lemma \ref{squarefinite}, there are $( q - 3 ) / 4$ nonzero values for 
      $( i + j )^2$ and $( q - 3 ) / 4$ nonzero values for $( i - j )^2$. But $i - j = 0$ also
      satisfies $1 - ( i - j )^2$ is nonzero square. Thus, we have $( q - 3 )
      / 2$ values for $i + j$ and $( q - 1 ) / 2$ values for $i - j$. From
      Lemma \ref{khac khong}, each pair $( i + j, i - j )$ gives us two points $Z$ and
      $Z'$. Thus, we have $( q - 1 ) ( q - 3 ) / 2$ points $Z$.      
      \item Suppose that $( i + j )^2 - 1$ and $1 - ( i - j )^2$ are nonsquare.  Similarly,
      we have $( q - 1 ) / 2$ values of $i + j$ and $( q - 3 ) / 2$ values of
      $i - j$. Thus, we also have $( q - 1 ) ( q - 3 ) / 2$ points $Z$.
    \end{enumerate}
  \end{enumerate}
  Therefore, the number of points $Z$ such that $Q ( X, Z ) = i^2$ and $Q ( Y, Z )
  = j^2$ is
  \[ 4 ( q - 2 ) + ( q - 1 ) ( q - 3 ) = q^2 - 5. \]
  But any point $Z$ is counted four times (for the four pairs $( i, j ), ( i, - j
  ), ( - i, j )$ and $( - i, - j )$), so $|\delta(X,Y)| = (q^2-5)/4$.
  
  b) Suppose that $( X, Y ) \in D_2$. Then $Q ( X, Y ) = -k^2$ for some $k\in F_q^{\ast}$.
  We count the number of points $Z$ in $\delta(X,Y)$.
  Let $Z$ be any point in $\delta ( X, Y )$ then there exists $i, j \in F_q^{\ast}$
  such that $Q ( X, Z ) = i^2$ and $Q ( Y, Z ) = j^2$. Without loss of generality, 
  we can divide both $i, j$ by $k$ and assume that $k = -1$. We have
  \begin{align*}
    4 f ( i^2, j^2, -1 ) &= 4 i^2 j^2 - ( -1 - i^2 - j^2 )^2 \\
    &= -( ( i + j )^2 + 1 ) ( 1 + ( i - j )^2 ) .
  \end{align*}
  There are three separate cases.
  \begin{enumerate}
    \item Suppose that $f ( i^2, j^2, -1 )$ is nonsquare. Then from Lemma                        
    \ref{khac khong}, there does not exist $Z$
    such that $Q ( X, Z ) = -i^2$ and $Q ( Y, Z ) = -j^2$.
    
    \item Suppose that $f ( i^2, j^2, -1 ) = 0$. This is impossible as $-1$ is not a square 
    number in $F_q$. 
    
    \item Suppose that $f ( i^2, j^2, -1 )$ is a nonzero square. Again, there are two subcases.
    \begin{enumerate}
      \item Suppose that $( i + j )^2 + 1$ is nonzero square and $1 + ( i - j )^2$ is nonsquare.
      From Lemma \ref{squarefinite}, there are $( q - 3 ) / 4$
      nonzero values for $( i + j )^2$ and $( q + 1 ) / 4$ nonzero values for 
      $( i - j )^2$. But $i + j = 0$ also
      satisfies $1 + ( i + j )^2$ is nonzero square. Thus, we have $( q - 1 )
      / 2$ values for $i + j$ and $( q + 1 ) / 2$ values for $i - j$. From
      Lemma \ref{khac khong}, each pair $( i + j, i - j )$ 
      gives us two points $Z$ and $Z'$. Thus, we have $( q - 1 ) ( q + 1 ) / 2$ points 
      $Z$.      
      \item Suppose that $( i + j )^2 + 1$ is nonsquare and $1 + ( i - j )^2$ is nonzero square.
      Similarly, we have $( q + 1 ) / 2$ values of $i + j$ and $( p - 1 ) / 2$ values 
      of $i - j$. Thus, we also have $( q - 1 ) ( q + 1 ) / 2$ points $Z$.
    \end{enumerate}
  \end{enumerate}
  Therefore, the number of points $Z$ such that $Q ( X, Z ) = i^2$ and $Q ( Y, Z )
  = j^2$ is $2(q-1)(q+1)/2=q^2 - 1$.
  But any point $Z$ is counted four times (for the four pairs $( i, j ), ( i, - j ), ( - i, j )$ and $( - i, - j )$) so $|\delta(X,Y)| = (q^2-1)/4$. This completes the proof of the theorem.
\end{proof}

From Theorem \ref{important}, we have $V_q$ is a strongly regular graph with parameters $\{(q^2-1)/2,(q^2-5)/4,(q^2-1)/4\}$; that is $V_q$ is $(q^2-1)/2$-regular, any two adjacent vertices have $(q^2-5)/4$ common neighbours and any two non-adjacent vertices have $(q^2-1)/4$ common neighbours. For any two vertices $a$ and $b$, there are precisely $(q^2-1)/4$ vertices $c \ne b$ joined to $a$ and not joined to $b$. We call $V_q$ \textbf{quadrance graph}.

\subsection{Maximal complete subgraphs}

Let $U$ be the vertex set of a subgraph of $V_q$ and $e(U)$ be the number of edges joining vertices in $U$. In Section 3.1, we see that the graph $V_q$ is a strongly regular graph with parameters $\{(q^2-1)/2,(q^2-5)/4,(q^2-1)/4\}$. The well-known Paley graph $P_{q^2}$ has the same parameter. Thus, we can follow Bollobas's proof for Paley graphs (see 
\cite{random graphs}, page 321-322) to show that
\begin{Theorem}
Let $U$ be a set of $u$ vertices of the Quadrance graph $V_q$. Then
\begin{equation*}
\left| e(U) - \frac{1}{2} \binom{u}{2} \right| \leqslant \frac{u}{4} \frac{q^2 - u}{q}.
\end{equation*}
\end{Theorem}
\begin{proof}
  We may assume that $0 < u < q^2$. Let $n ( U )$ be the number of unordered
  triples $( X, Y, Z )$ with $X, Y \in U$ and $( X, Z ), ( Y, Z ) \in E ( V_q
  )$. By Theorem 4, for each pair $( X, Y )$ of adjacent vertices of $U$ there
  are $( q^2 - 5 ) / 4$ such triples and for each pair $( X, Y )$ of
  non-adjacent vertices of $U$ there are $( q^2 - 1 ) / 4$ such triples. Thus,
  we have
  \[ n ( U ) = e ( U ) \left( \frac{q^2 - 5}{4} \right) + \left\{ \binom{u}{2} - e
     ( U ) \right\} \left( \frac{q^2 - 1}{4} \right) = \binom{u}{2} \left(
     \frac{q^2 - 1}{4} \right) - e ( U ) . \]
  Let $d = 2 e ( U ) / u$ be the average degree in the induced subgraph $G [ U
  ]$. There are $[ \{ ( q^2 - 1 ) / 2 \} - d ]$ edges joining $U$ to $V - U$.
  So on average, a vertex of $U$ is joined to $d$ vertices in $U$, and a
  vertex in $V - U$ is joined to $d_0 = \{ u / ( q^2 - u ) \} [ \{ ( q^2 - 1 )
  / 2 ) \} - d ]$ vertices in $U$. Thus, we have
  \[ n ( U ) \geqslant u \binom{d}{2} + ( q^2 - u ) \binom{d_0}{2}. \]
  This implies
  \[ \frac{u d ( d - 1 )}{2} + \frac{u}{2} \left( \frac{q^2 - 1}{2} - d
     \right) \left\{ \frac{u}{q^2 - u} \left( \frac{q^2 - 1}{2} - d \right) -
     1 \right\} \leqslant \frac{u ( u - 1 ) ( q^2 - 1 )}{8} - \frac{d u}{2} .
  \]
  Multiplying by $u ( q^2 - u ) / 2 q^2$ and rearranging, we have
  \[ \left\{ \frac{d u}{2} - \frac{1}{2} \binom{u}{2} \right\}^2 \leqslant
     \frac{u^2}{16 q^2} ( q^2 - u )^2 . \]
  Hence
   \[ \left| e(U) - \frac{1}{2} \binom{u}{2} \right| \leqslant \frac{u}{4} \frac{q^2 - u}{q}.  \]
This completes the proof.
\end{proof}

If $U$ spans a complete subgraphs or an empty subgraphs of $V_q$ and $|U|=u$ then 
\[\left| e(U) - \frac{1}{2} \binom{u}{2} \right| = \frac{1}{2} \binom{u}{2} \leqslant \frac{u}{4} \frac{q^2 - u}{q}.\]
Hence $u \leqslant q$. It is easy to see that the equality holds when $U$ is a set of $q$ points in a line. Inspections for $q= 3, 7$ support for the converse statement. We propose the following conjecture.
\begin{Conjecture}
If $U$ spans a complete subgraphs or an empty subgraphs of order $q$ of $V_q$ then $U$ is a line in $F_q \times F_q$.
\end{Conjecture}

\section{Companion Results} \label{companion results}

In this section, we suppose that $q \equiv 1 \:( \tmop{mod} 4 )$. We present companion results of previous sections in this case. Some proofs of the results in this section require extra works but they are basically similar to proofs in previous sections. Therefore, we obmit the details. The following lemma is an analogue to Lemma 1.

\begin{Lemma} \label{comp 1} Let $q \equiv 1 \:( \tmop{mod} 4 )$ be an odd prime power. Then
$A_q = (q-5)/5$, $B_q = (q-1)/4$ and $C_q = D_q = (q-1)/4$.
\end{Lemma}

Using Lemma \ref{comp 1}, we have the following result, which is an analogue to Lemma 2. 

\begin{Lemma} \label{comp 2}
  Let $C$ be a circle with nonzero quadrance in $F_q^2$. Then $C$ has exactly $N = q-1$ points.
\end{Lemma}

Lemma 4 still holds in this case. From Lemma \ref{number points} and Lemma \ref{comp 2}, it is clear that the circle $C_0(X)$ is a single point when $q = 4 l + 3$ for some integer $l$. But for $q = 4l+1$, the circle $C_0(X)$ consists of $2q-1$ points. This property causes some troubles in calculations (see \cite{hypergroup}) which can be resolved by redefine the circle $C_0(X)$. We redefine $C_0(X) = \{X\}$ and $C_q(X) = \{ Y \in F_q \times F_q \;|\; Q(Y,X)=0, Y \ne X\}$. We have the following lemmas.

\begin{Lemma} \label{bang khong}
For any $i, j \ne 0$ in $F_q$. Suppose that $X, Y$ are two distinct points in $F_q^2$ with $Q(X,Y)= 0$. Then the circle $C_i(X)$ intersects $C_j(Y)$ if and only if $i \ne j$. Furthermore, if $i \ne j$ then two circles intersect at only one point.
\end{Lemma}

\begin{Lemma} \label{bang nhau}
  Let $i \in F_q^{\ast}$ then there exists $j \in F_q^{\ast}$ such that $f (
  i, i, j )$ and $f ( j, j, i )$ are square in $F_q$.
\end{Lemma}

\begin{proof}
  There are three cases.
  \begin{enumerate}
    \item Suppose that $3$ is square in $F_q$. Then $4 f ( i, i, i ) = 4 i^2 -
    ( i - i - i )^2 = 3 i^2$ is square in $F_q$. Thus, we can choose $j = i$.
    
    \item Suppose that $5$ is square in $F_q$. We have
    \[ 4 f ( i, i, - i ) = 4 f ( - i, - i, i ) = - 5 i^2 . \]
    Since $- 1$ is square in $F_q$ so we can choose $j = - i$.
    
    \item Suppose that neither $3$ nor $5$ is square in $F_q$. Then $15$ is a
    square in $F_q$. We have
    \begin{align*}
      4 f ( i, i, - 2 i ) &= 4 i^2 - ( - 6 i - i - i )^2 = - 60 i^2,\\
      4 f ( - 6 i, - 6 i, i ) &= 144 i^2 - ( i + 6 i + 6 i )^2 = - 25 i^2 .
    \end{align*}
    Thus, both $f ( i, i, - 2 i )$ and $f ( - 2 i, - 2 i, i )$ are square in
    $F_q$. Hence we can choose $j = - 6 i$.
    
  \end{enumerate}
  
  This concludes the proof.
\end{proof}

Lemmas \ref{khac khong}, \ref{function}, \ref{bang khong} and \ref{bang nhau} can be used to obtain an analogue result for Theorem \ref{polygon}.
\begin{Theorem}  
  If $q \equiv 1 ( \tmop{mod} 4 )$ then for any $n \geqslant 5$ and $a_1,
  \ldots, a_n \in F_q^{\ast}$, there exists a polygon $A_1 \ldots A_n$ in
  $F_q^2$ such that $Q ( A_i, A_{i + 1} ) = a_i$ for all $i \in \{ 1, \ldots,
  n \}$. 
\end{Theorem}

\begin{proof}
	Note that the condition $n \geqslant 5$ is necessary in this case. A
  simple illustration is that there does not exist a quadrangle with quadrances
  $1, 1, 1$ and $3$ in $F_5^2$. Now we prove that for any $n \geqslant 5$ and $a_1,
  \ldots, a_n \in F_q^{\ast}$ then there exists a polygon $A_1 \ldots A_n$
  in $F_q^2$ such that $Q ( A_i, A_{i + 1} ) = a_i$ for all $i$. It is clear
  that we only need to prove the result for $n = 5$. The case $q = 5$ is easy to check,
  so we assume that $q > 5$ (then $q \geqslant 3^2 = 9$). We consider two cases.
  
  Case 1: Suppose that $a_1 = a_2 = \ldots = a_5$. From Lemma \ref{bang nhau}, there exists $j \in F_q^{\ast}$ such that
  $f(a_1,a_1,j)$ and $f(j,j,a_1)$ are square in $F_q$. From Lemma \ref{khac khong}, there exists a triangle $A_1A_2A_4$ 
  with $Q(A_1A_2)=a_1$ and  $Q(A_2,A_4)=Q(A_1,A_4) = j$. We have $Q(A_2,A_4)=j$ so from Lemma \ref{khac khong}, there
  exists a point $A_3$ such that $Q(A_2A_3)=Q(A_3A_4)=a_1$. Similarly, there exists $A_5$ such that 
  $Q(A_4A_5)=Q(A_1A_5)=a_1$. And we have a pentagon $A_1A_2A_3A_4A_5$ as required.
  
  Case 2: Suppose that $a_1 \neq a_2$. Since Lemma \ref{function} still holds,
  there are at least $( p - 1 ) / 2 >
  3$ values of $k$ such that $f ( a_3, a_4, k )$ is square. Thus, we can
  choose some $k \neq 0, a_5$. Now we have $k \neq a_5$ and $a_1 \neq a_2$. Hence
  by Lemma \ref{bang khong}, for any $X \ne Z \in F_q^2$ such that $Q ( X, Z ) = 0$ then
  there exists $Y, V$ such that $Q ( X, Y ) = a_1$, $Q ( Y, Z ) = a_2$, $Q (
  Z, V ) = k$ and $Q ( V, X ) = a_5$. But $f ( a_3, a_4, k )$ is
  square so by Lemma \ref{khac khong}, there exists $T$ such that $Q ( Z, T ) = a_3$ and
  $Q ( T, V ) = a_4$. Hence we can choose the pentagon $A_1 \ldots A_5 \equiv X
  Y Z T V$. This concludes the proof.
\end{proof}

We have an open question.

\begin{Problem}
For $q \equiv 1 \:( \tmop{mod} 4 )$ is an odd prime power. Find the necessary and sufficient condidtion for $i, j,k, l \in F_q$ such that there exist a quadrangle with quadrances $i, j, k$ and $l$.
\end{Problem}

We can define quadrance association scheme on $F_q^2$ similarly as in Section 2. But note that in this case, we need to partition $F_q^2$ into $q+1$ partitions. Precisely, we have the following theorem. 

\begin{Theorem}
  Let $\Omega = F_q^2$, we partition $\Omega \times \Omega$ into $q+1$ subsets $C_0, C_1, \ldots, C_q$ such that for $X, Y$ in $F_q^2$, $(X, Y ) \in C_i$ if and only if $X \in C_i(Y)$. This partition gives us an association scheme on $F_q^2$. Furthermore, for $i, j, k \in F_q$, let $f ( i, j, k ) = i j - ( i - j - k )^2 / 4$ then the
  intersection numbers of the scheme are
  \begin{equation*}
  p_{ij}^0 =
   \begin{cases}
    0  & \text{if}\; \; i \neq j,\\
    1  & \text{if}\; \; i = j = 0,\\
    2(q-1)  & \text{if}\; \; i = j = q,  \\ 
    q-1 & \text{otherwise},
   \end{cases}
  \end{equation*}
    \begin{equation*}
  p_{ij}^q =
   \begin{cases}
    1  & \text{if}\; \; i \neq j,\\
    q-2  & \text{if}\; \; i = j = q,\\
    0 & \text{otherwise},
   \end{cases}
  \end{equation*}
and
  \begin{equation*}
  p_{ij}^k =
   \begin{cases}
    0  & \text{if}\; \; f(i,j,k)\; \text{is non-square},\\
    1  & \text{if}\; \; f(i,j,k) = 0,\\
    2 &\text{if}\; \; f(i,j,k)\; \text{is nonzero square}.
   \end{cases}
  \end{equation*}
for $k \neq 0, q$.

\end{Theorem}

We neither know an analogue to Theorem \ref{induced scheme} for the association scheme above nor any strongly regular quadrance graph in this case.

\end{document}